\documentclass[11pt,draftcls,onecolumn]{ieeeconf}
\usepackage{graphicx} 
\usepackage{amsmath} 
\usepackage{amssymb}  
\usepackage{amsfonts}
\usepackage{color}
\usepackage[english]{babel}
\usepackage[latin1]{inputenc}

\newcommand{\G}{\mathcal{G}}

\newcommand{\xmin}{x^{\star}}
\newcommand{\xtrue}{\alpha}

\newcommand{\dy}{\delta_y}
\newcommand{\dA}{\delta_A}
\newcommand{\DA}{\Delta_A}
\newcommand{\dAj}{\delta_{A,j}}
\newcommand{\dAi}{\delta_{A,i}}
\newcommand{\Dy}{\Delta_y}
\newcommand{\qA}{\mathcal{Q}(A)}
\newcommand{\qy}{\mathcal{Q}(y)}

\newcommand{\can}{\beta} 
\newcommand{\FS}{\mathcal{D}}
\newcommand{\su}{\mathcal{S}}

\newcommand{\R}{\mathbb{R}}

\newcommand{\one}{\mathbf{1}}

\newtheorem{theorem}{Theorem}

\newtheorem{assumption}{Assumption}

\usepackage{graphics} 
\usepackage{epsfig} 

\usepackage{mathtools}

\usepackage{tikz}
\usetikzlibrary{shapes,arrows,calc,positioning}

\tikzset{
	block/.style    = {draw, thick, rectangle, minimum height = 2.5em,
		minimum width = 2.5em, node distance = 1.7cm},
     sum/.style      = {draw, circle, radius = 0.5cm, node distance = 1.3cm}, 
     input/.style    = {coordinate}, 
     output/.style   = {coordinate} 
 }

\newtheorem{remark}{\noindent \textbf{Remark}}{\normalfont}{\normalfont}
\newtheorem{result}{\noindent \textbf{Result}}{\normalfont}{\normalfont}
\newtheorem{problem}{\noindent \textbf{Problem}}{\normalfont}{\normalfont}
 {\normalfont }{\normalfont}
{\normalfont}{\normalfont}

\begin{document}
\title{A linear programming approach to \\ sparse linear regression with quantized data}
\author{V. Cerone, S. M. Fosson$^{*}$, D. Regruto \thanks{$^{*}$ Corresponding author. The authors are with the Dipartimento di Automatica e Informatica,
Politecnico di Torino,
    corso Duca degli Abruzzi 24, 10129 Torino, Italy;
    e-mail: vito.cerone@polito.it, sophie.fosson@polito.it, diego.regruto@polito.it;
}
}

\maketitle
\begin{abstract}The sparse linear regression problem is difficult to handle with usual sparse optimization models when both predictors and measurements are either quantized or represented in low-precision, due to non-convexity. In this paper, we provide a novel linear programming approach, which is effective to tackle this problem. In particular, we prove theoretical guarantees of robustness, and we present numerical results that show improved performance with respect to the state-of-the-art methods.
\end{abstract}

\section{Introduction}\label{sec:intro}
Sparse optimization refers to those optimization problems where the solution is encouraged to be sparse, i.e., to have few non-zero components. This research area has dramatically increased in the last decades in many different fields including signal processing, machine learning, and system identification. In signal processing, the widespread presence of signals that admit sparse representations has promoted the research on new sparse optimization problems and algorithms, see, e.g., \cite{don06,fou13}. In machine learning and system identification, sparsity is desirable to reduce as much as possible the complexity of the estimated models. In the literature, sparsity is exploited for the identification of linear systems, see, e.g., \cite{gu09,tot11,tothsysid12}; non-linear functions in  \cite{novacc11}; polynomial models in \cite{cal14ifac}; time-varying  systems in \cite{san11,fox18cdc}. 

A popular paradigm in sparse optimization is given by the case where the collected measurements $y$ are linearly related, through a predictor matrix $A$, to the sparse signal or vector of parameters $x$ to be estimated. This paradigm leads to the problem of computing the sparsest solution of the linear system of equation $Ax=y$, $x\in\R^n$, $y\in\R^m$, $A\in\R^{m,n}$. In the literature, this  problem is generically referred to as sparse linear regression.  The undetermined case, $m<n$, has attracted the attention of many researchers, whose work originated the theory of compressed sensing (CS) \cite{don06,fou13}. Finding the sparsest solution of an undetermined linear system is an NP-hard problem. However, it is well known that sparsity can also be achieved by minimizing the $\ell_1$-norm of $x,$ under suitable constraints accounting for the linear structure of the problem, which makes the problem convex. Specifically, the minimization of the $\ell_1$-norm of $x$ subject to $Ax=y$ is known as Basis Pursuit  \cite[Chapter 4]{fou13}. When a measurement noise is present, the constraint is generally formulated as $\|Ax-y\|_p\leq \epsilon$, where $\|\cdot\|_p$ is a suitable norm and $\epsilon>0$ is a known bound; this is referred to as Basis Pursuit Denoising (BPDN$_p$).

In the literature, $p=2$ and $p=\infty$ are the most common choices. In particular, BPDN$_2$ is very popular, as it is suitable to cope with Gaussian noise, which is the typical model in a number of applications, such as  transmission systems. The choice $p=2$ provides solutions that are more tolerant to possible outliers, since it bounds the mean energy of the error.

The case $p=\infty$  was first analyzed in \cite{don06infty}, where results on robustness to noise are proven based on the coherence properties of $A$, and has been recently retrieved to deal with quantized or low-precision measurements in the CS setting. When $y$ is quantized, in fact, there is a bounded error on each component $y_i$, $i=1,\dots,m$, which makes  the $\ell_{\infty}$-norm description more suitable than the $\ell_2$ one, as illustrated in \cite{las11,val15}. In particular, the $\ell_{\infty}$-norm supports the \emph{consistency} principle: the measurements obtained from the recovered signal lie in the same quantization intervals of the observed measurements, as considered in \cite{1bit,jac11,las11,val15}.

As to CS, the study of quantization is strongly motivated from the practical point of view. As a matter of fact, the CS paradigm moves the computational burden from the acquisition-compression phase (which simply consists in computing $y=Ax$) to the recovery phase. This feature is successfully exploited in systems where signals' acquisition is performed by remote devices with reduced computational capability, e.g., either space probes or environmental sensors, while recovery is performed in powerful computational centers. However, for transmission purposes, measurements are often not only compressed, but also quantized. We refer the interested reader to \cite{bou15,shi16} for a complete overview on CS with quantized measurements.

In many applications, also $A$ is generated by the remote device and  has to be  transmitted. Therefore, it is more realistic to assume that also $A$ undergoes quantization. It is worth noticing that, in some cases, $A$ can be designed on purpose, then it may be quantized in its original form, which prevents the addition of a quantization error. For example, Bernoulli matrices (whose entries are binary) are suitable for CS. However, in most of the applications, $A$ is not arbitrarily chosen, instead it depends either on the hardware of the remote device or on the physical nature of the problem.

The quantization of $A$ plays the role of an undesired perturbation in the recovery phase. Inaccuracies in $A$ are a serious drawback because they lead to an uncertain sparse linear regression problem which is not convex anymore. The problem of perturbed $A$ in sparse optimization and CS is addressed in  \cite{her10,zhu11,yan12}. 

The goal of this paper is to tackle the problem of sparse linear regression when both $A$ and $y$ are quantized, with main focus on the CS setting. To the best of our knowledge, this joint problem has been considered only in \cite{gur18}, while, as already mentioned, the two single problems of quantized $y$ and perturbed $A$  have already been studied. In \cite{gur18}, the normalized iterative hard thresholding algorithm \cite{blu10} is adapted to tackle the  problem of quantized $A$ and $y$. The algorithm is tested in the case of a specific stochastic quantization function, with possible applications in the radio astronomy framework. In this setting, a robustness result is proved, which shows that the mean recovery error is controlled by the quantization error.

In this paper, we propose a novel approach to this problem, which can be applied in the presence of any quantization function with bounded error. We extend the formulation of BPDN$_{\infty}$ and the results in \cite{don06infty} to the case where also $A$ is quantized. First we recast the considered problem  into the framework of the static error-in-variables  estimation considered in \cite{cer93}. Then, by exploiting results from \cite{cer93}, we show that the solution can be computed by solving a suitable number of linear programming (LP) problems. 
The paper is organized as follows. In Section \ref{sec:ps}, we formally state the problem and specify the considered assumptions. In Section \ref{sec:lp}, we introduce the novel LP  formulation.  In Section \ref{sec:ar}, we prove theoretical results on the robustness of the proposed approach. In Section \ref{sec:nr}, we show some numerical simulations that support the efficiency of the proposed method with respect to the state-of-the-art. Finally, some conclusions are drawn in Section \ref{sec:con}.
\section{Problem statement}\label{sec:ps}
Let us consider a device performing compressed data acquisition according to the equation $Ax=y$, $A\in\R^{m,n}$, where $m<n$. 
The device is assumed to transmit  quantized/low-precision versions of $A$ and $y$ to a recovery center. We denote by $\qA$ and $\qy$ the quantized versions of $A$ and $y$, respectively.  The problem considered in this work is to recover a sparse $x$ such that $Ax=y$, given $\qA$ and $\qy$. The quantization strategy is supposed to be unknown, though a bound on the maximum quantization error is given. Apart from the quantization, for simplicity, no other sources  of uncertainty/noise are considered here, although
such extension is under investigation. More precisely, we formulate the following optimization problem.\\
\\
\begin{problem}
Given $\qy$, $\qA$, $\DA>0$, and $\Dy>0$,
\begin{equation}\label{theproblem}
 \begin{split}
  \min_{x\in\R^n}\|x\|_1~~~ \text{s. t. }& y=A x\\
  &\qy=y+\dy\\
  &\qA=A+\dA\\
  &\|\dy\|_{\infty}\leq \Dy\\
  &\|\dA\|_{\infty}\leq \DA.
 \end{split}
\end{equation}\\
\end{problem}

In the rest of the paper, we assume that the signs of the components of $x$ are known, in the sense that for each $x_i$, $i=1,\dots,n$, we know if it is either non-negative or non-positive. More specifically, without loss of generality, we assume the non-negativity of $x$. 
\begin{assumption}\label{ass1}
$x_i\geq 0$ for all $i=1,\dots,n$.\\
\end{assumption}
This assumption naturally occurs in a number of applications, such as localization problems \cite{bay15}, image processing \cite{fra17}, and power allocation \cite{cal16}. However, extensions to more general classes are possible and will be studied in future work.
Here, we only observe that if the signs are not known, one can split the problem in $2^n$  LP problems, trying all the possible combinations of signs. Therefore, there is a way to compute the global minimum of the general non-convex problem, though computationally not efficient for large $n$. 
In future work, we will also analyze possible pre-processing methods to obtain prior  information on the signs from available data.

\section{A linear programming approach}\label{sec:lp}
Thanks to the following result, we show that, under Assumption \ref{ass1}, Problem 1 is equivalent to an LP problem. \\
In the following, given two vectors $a,b\in\R^n$, we write $a\succeq b$  to indicate that $a_i \geq b_i$ for each $i=1,\dots,n$. We denote by $I_n\in\R^{n,n}$ the identity matrix. Moreover, $\one_n:=(1,1,\dots,1)^T\in\R^n$, where $T$ is for transpose. 
\\
\begin{result}
Under Assumption \ref{ass1}, Problem 1 can be equivalently formulated  as the following LP problem:
\begin{equation}\label{sys3}
 \begin{split}
  \min_{x\in\R^n}\|x\|_1~&\text{ s. t. } C x \preceq c\\
  &\text{ where }\\
  &C= \left(\begin{array}{c}
         \qA-\DA\one_m\one_n^T\\
         -\qA-\DA\one_m\one_n^T\\
        \end{array}\right)\in\R^{2m,n}\\
    &c= \left(\begin{array}{c}
         \qy+\Dy\one_m\\
         -\qy+\Dy\one_m\\
        \end{array}\right)\in\R^{2m}.\\
        \end{split}
\end{equation}
\\
\end{result}
Proof of Result 1 is obtained by first noticing that Problem 1 can be equivalently rewritten in a more compact way as follows:
\begin{equation}\label{sys2}
 \begin{split}
  \min_{x\in\R^n}\|x\|_1~\text{ s. t. } &\|\qy-(\qA-\dA)x\|_{\infty}\leq \Dy\\
  &\|\dA\|_{\infty}\leq \DA.
 \end{split}
\end{equation}
Model \eqref{sys3} is obtained by applying the results about bounded errors-in-variables identification of static linear systems presented in \cite{cer93} to the set of constraints of problem \eqref{sys2}, under Assumption \ref{ass1}. Model \eqref{sys3} with Assumption \ref{ass1} is an LP problem whose solution is straightforward. Moreover, its formulation is compliant with the quantization consistency principle.
\\
\begin{remark}
It is worth noting that, by suitably exploiting results in \cite{cer93}, Result 1 can be extended to cover the general case where the signs of $x\in\R^n$ are unknown, although in that case the solution is obtained by solving a larger number of LP's. However, this general case is outside the scope of this conference contribution and will be the subject of a future work.
\\
\end{remark}
\begin{remark}
We notice that a working hypothesis similar to Assumption \ref{ass1} is considered in \cite[Theorem 6]{yan12}. In \cite{yan12}, the perturbation on $A$ is assumed to have a specific structure, namely, it is equal to $B\text{diag}(\beta_0)$, where  $B\in\R^{m,n}$ is known, and $\beta_0\in[-r,r]^n$, $r>0$, is unknown. In other terms, the direction of each column of the perturbation is known, and the dimension of the unknown is reduced from $mn$ to $n$. The proposed model \cite[Equation 11]{yan12} is a BPDN$_2$ with an additive perturbation bounded in the $\ell_{\infty}$-norm. This model is biconvex and can be approached with alternating minimization, which only achieves a local minimum. However, in \cite[Theorem 6]{yan12}, it is observed that  if $x_i\geq 0$, for all $i=1,\dots,n$, the problem admits convex formulation, which guarantees to get the global minimum.
\\
\end{remark}

\section{Analysis of robustness}\label{sec:ar}
In this section, we show that Model \eqref{sys3} with Assumption \ref{ass1} is robust, \emph{i.e.}, the distance between its solution and the original signal is bounded by a quantity $T>0$, which is controlled by $\DA$ and $\Dy$. In particular, if $\DA=0$, the robustness result in \cite{don06infty} is obtained.


The result that we now prove is based on the \emph{mutual coherence} of $\qA$, which is defined as:
\begin{equation}
\mu:= \max_{j\neq h} |\qA_j^T \qA_h|~~~~j,h\in\{1,\dots,n\} 
\end{equation}
where the index $j$ denotes the $j$-th column. 

From the theory on underdetermined linear systems and CS, it is well known that a sufficiently small coherence  guarantees a successful sparse linear regression, up to errors due to noise \cite{fuc04,fuc05,don06infty}. The following theorem shows that sparse linear regression from quantized data is successful when coherence is sufficiently small, up to quantization errors.
\begin{theorem}\label{theo:robust}
Let $\qy=A \alpha+\dy$, where the unknown $\alpha\in\R^n$ has $k\ll n$ non-zero components. $\qA=A+\dA$ and $\qy$ are known. 
Let us assume that: $\|\qA_j\|_2\leq \rho$ for some $\rho>0$, for any $j=1,\dots,n$, $2(\DA+\rho)^2<2-\mu-\rho^2$,  and
$$k\leq\frac{1}{2} \frac{2-\rho^2+\mu}{ \mu+\DA^2+2\DA\rho  + (\rho\sqrt{m}+\DA m)2\Dy/T}.$$
Then, the solution $\xmin$ of problem \eqref{sys3} is robust, that is,
$$\|\xmin-\xtrue\|_1<T.$$
\end{theorem}
\begin{proof}
Let 
$$\FS:=\{x\in\R^n: \|\qy-A x\|_{\infty}\leq \Dy\}$$
be the feasible set of Model \eqref{sys3}. By definition, $\alpha\in\FS$.
Let us consider the subset  $\G:= \{\can\in \FS: \|\can-\xtrue\|_1\geq T\}$. We then prove that for any $\can\in\G$ we have $\|\can\|_1 \geq \|\xtrue\|_1$, which means that there is no solution in $\G$. This implies that all the solutions are in $\FS\setminus \G$, which proves the thesis.
Thus, we study the problem:
\begin{equation}\label{aux1}
\min_{\can} \|\can\|_1-\|\xtrue\|_1 \text{ s. t. } \can\in\G,~\xtrue\in\FS
\end{equation}

Let  $w:=\can-\xtrue.$
As illustrated in \cite{don01}, it is straightforward to prove that
\begin{equation}\label{use_w}
\|\can\|_1-\|\xtrue\|_1\geq\|w\|_1-2\sum_{h\in \su}|w_h|
\end{equation}
where $\su$ is the support of $\alpha$.
Since $\xtrue,\can \in \FS$, we  obtain 
\begin{equation}\label{contraint1}
\|A w\|_{\infty}\leq 2\Dy
\end{equation}
from which we also get $\|A w\|_{2}\leq 2\Dy\sqrt{m}.$
By assuming $\|\qA_j\|_2\leq \rho$ for any column $j$, we have:
\begin{equation}\label{eq:useful_inequality}
\begin{split}
|A_j^T A w|&\leq \|A_j\|_2 \|A w\|_{2} \leq (\rho+\DA\sqrt{m}) (2\Dy\sqrt{m}).
\end{split}
\end{equation}
where we use the triangle inequality $\|A_j\|_2\leq \|\qA_j\|_2+\|\dAj\|_2$.
Now, we notice that
\begin{equation}\label{eq:splittino}
\begin{split}
|w|&= |(A^T A -A^T A+I_{n}) w|  \preceq |A^T A w|+ |A^T A-I_{n}| |w|.
\end{split}
\end{equation}
where $|w|=(|w_1|\,\dots,|w_n|)^T$. From \eqref{eq:useful_inequality}, we obtain a bound for $|A^T A w|$. 
Furthermore, for the off-diagonal elements of $A^T A$, we have:
\begin{equation*}
\begin{split}
A_i^T A_j &= (\qA_i^T+\dAi^T)(\qA_j+\dAj)\\ 
&\leq \mu+\DA^2+2\DA\rho.
\end{split}
\end{equation*}
Similarly, for the diagonal elements, we obtain: $$A_i^T A_i\leq \rho^2+\DA^2+2\DA\rho.$$
Therefore,
\begin{equation*}
\begin{split}
|A^T A-I_{n}|&\preceq (\mu+\DA^2+2\DA\rho)(\one_{n,n}-I_n)\\
&~~+(\rho^2+\DA^2+2\DA\rho)I_n-I_n\\
&=(\mu+\DA^2+2\DA\rho)\one_{n,n}+(\rho^2-\mu-1)I_n
\end{split}
\end{equation*}
where $\one_{n,n}:=\one_n\one_n^T$.
Coming back to \eqref{eq:splittino},
\begin{equation}\label{eq:splittino_bound}
\begin{split}
|w|&\preceq  (\rho+\DA\sqrt{m}) (2\Dy\sqrt{m})\one \\
&+(\mu+\DA^2+2\DA\rho)\one_{n,n}|w| +(\rho^2-\mu-1)|w|.
\end{split}
\end{equation}
By assuming $2-\rho^2+\mu\geq 0$, we can write
\begin{equation*}
|w|\preceq\frac{(\rho\sqrt{m}+\DA m)2\Dy}{2-\rho^2+\mu}\one  + \frac{\mu+\DA^2+2\DA\rho}{2-\rho^2+\mu}\one_{n,n}|w|.
\end{equation*}
Finally,
\begin{equation*}
\left(\hspace{-0.12cm}I_n\hspace{-0.12cm}-\frac{\mu+\DA^2+2\DA\rho}{2-\rho^2+\mu}\one_{n,n}\hspace{-0.1cm}\right)|w|\preceq\frac{(\rho\sqrt{m}+\DA m)2\Dy}{2-\rho^2+\mu}\one_n.
\end{equation*}
Let $v=|w|.$ We can now rewrite problem \eqref{aux1}, using \eqref{use_w}, as follows:
\begin{equation}\label{aux2}
\begin{split}
&\min_{v\in\R^n} (\one_n-2\one_{n}^{\su})^T v \text { s. t. }\\
&\left(\hspace{-0.1cm}I_n\hspace{-0.1cm}-\frac{\mu+\DA^2+2\DA\rho}{2-\rho^2+\mu}\one_{n,n}\hspace{-0.1cm}\right)v\preceq\frac{(\rho\sqrt{m}+\DA m) 2\Dy}{2-\rho^2+\mu}\one_n\\
&\one_n^T v \geq T,~v\succeq 0
\end{split}
\end{equation}
where $\one_{n}^{\su}$ is the $n$-dimensional column vector with entries equal to 1 in the positions of the support $S$ of $\alpha$, and 0 otherwise.

As in \cite[equations (20)-(21)]{don06infty}, we consider the dual problem:
\begin{equation}\label{aux2dual}
\begin{split}
&\max_{u\in\R^n} -\frac{(\rho\sqrt{m}+\DA m) 2\Dy}{2-\rho^2+\mu}\one_n^T u+ Tu_0~~\text { s. t }\\
&\one_n u_0-\left(\hspace{-0.1cm}I_n\hspace{-0.1cm}-\frac{\mu+\DA^2+2\DA\rho}{2-\rho^2+\mu}\one_{n,n}\hspace{-0.1cm}\right)u\preceq \one_n-2\one_{n}^{\su}\\
&u\succeq 0, u_0 \geq 0\\
\end{split}
\end{equation}
Exploiting the zero duality gap between primal and dual in LP problems \cite{lue15}, if \eqref{aux2dual} has solution that originates a positive penalty, the penalty is positive also for \eqref{aux2}, which is our final aim. 
From this point, the thesis can be obtained following the same procedure used in \cite[pages 518-519]{don06infty}, since problem \eqref{aux2dual} is analogous to problem (21) in \cite{don06infty}, with different constants. We omit the details for brevity. We just notice that the constraint $2(\DA+\rho)^2<2-\mu-\rho^2$ is necessary to fulfill equations (22)-(23) in \cite{don06infty}. 
\end{proof}

We remark that if $\DA=0$, Theorem \ref{theo:robust} provides the same bound of Theorem 3 in \cite{don06infty} (with $\delta=\epsilon$).

It is worth noticing that, in CS, coherence-based analyses \cite{fuc04,fuc05}  have  the drawback of providing less tight bounds with respect to other properties, such as the restricted isometry property (RIP) \cite{fou13}. Nevertheless, RIP is difficult to assess for a specific matrix (in the literature, RIP is proved for some classes of random matrices). Coherence, instead, can be easily computed for any matrix. This work provides results in terms of coherence, while future extensions might envisage other properties.
\section{Numerical results}\label{sec:nr}
\begin{figure*}
\centering
\includegraphics[width=0.47\textwidth]{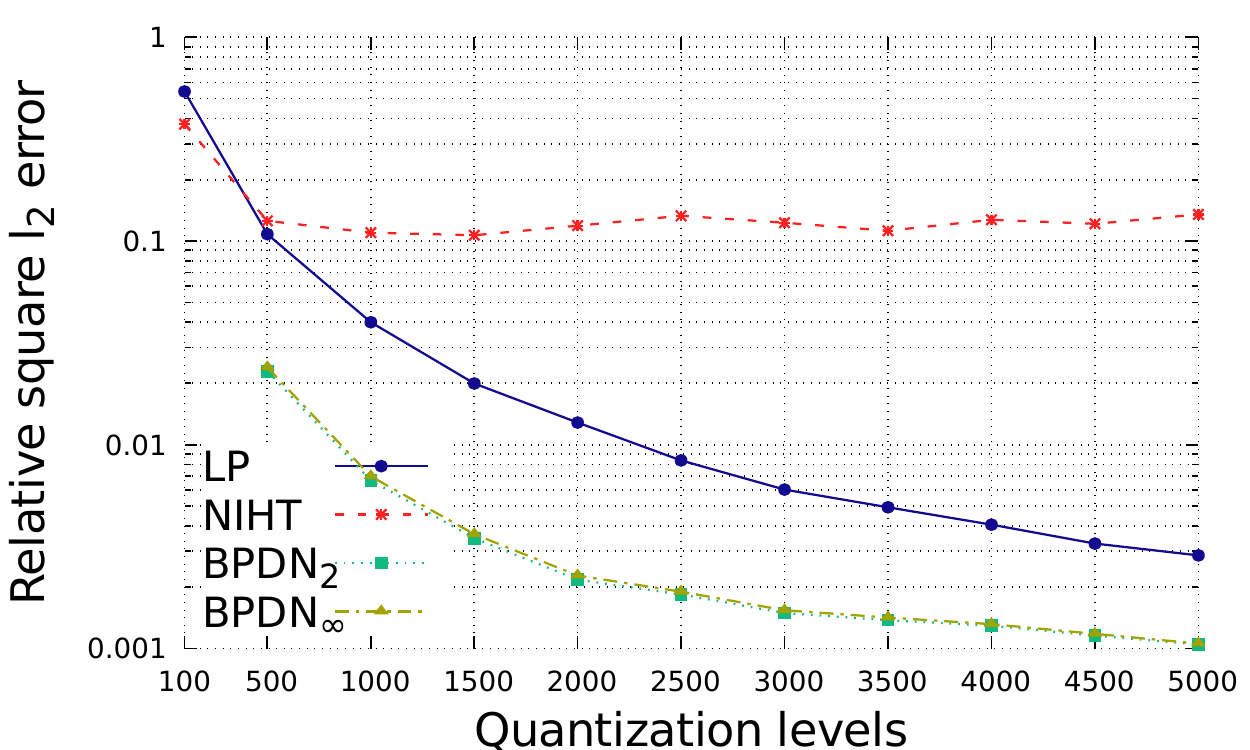}\quad
\includegraphics[width=0.47\textwidth]{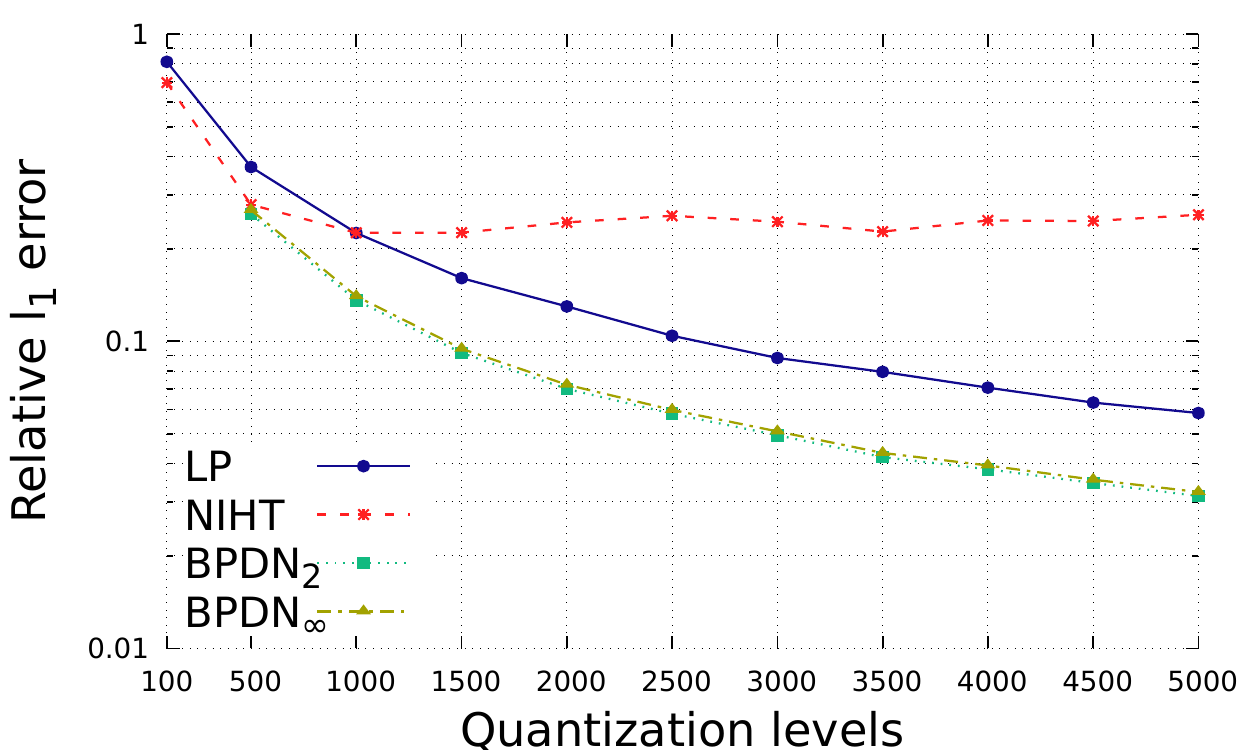}
\includegraphics[width=0.47\textwidth]{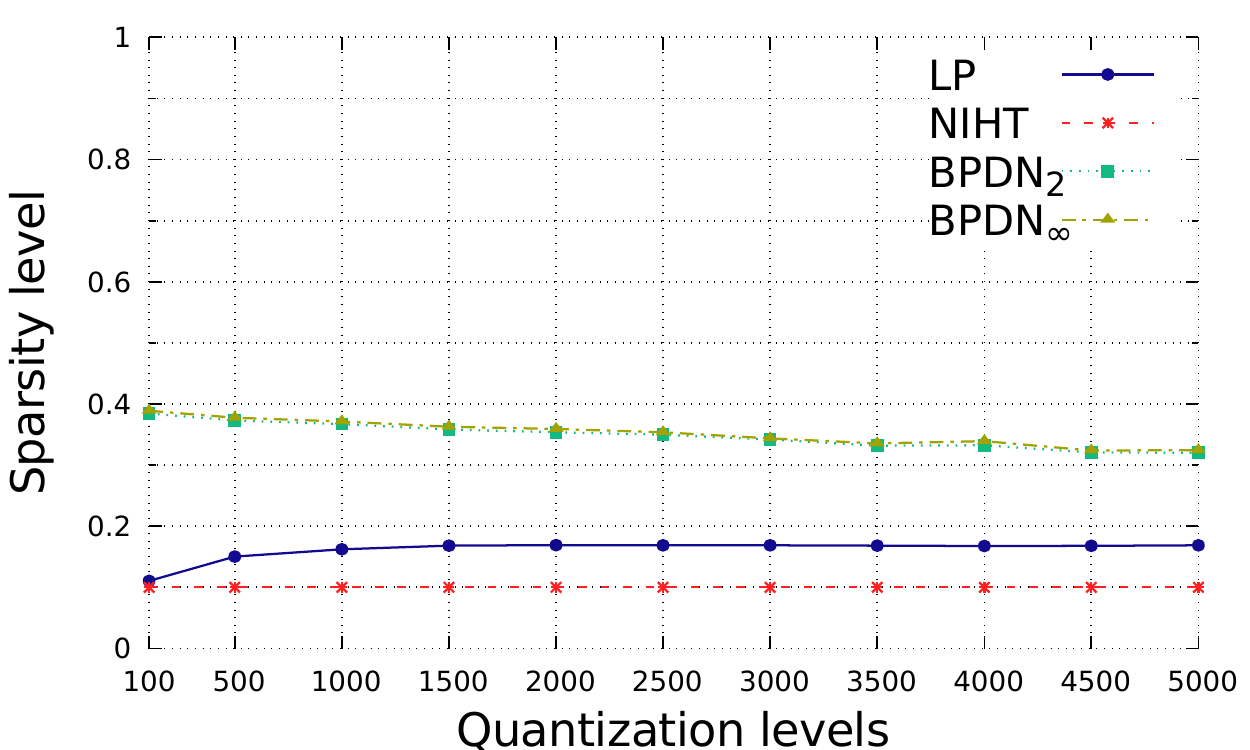} 
\includegraphics[width=0.47\textwidth]{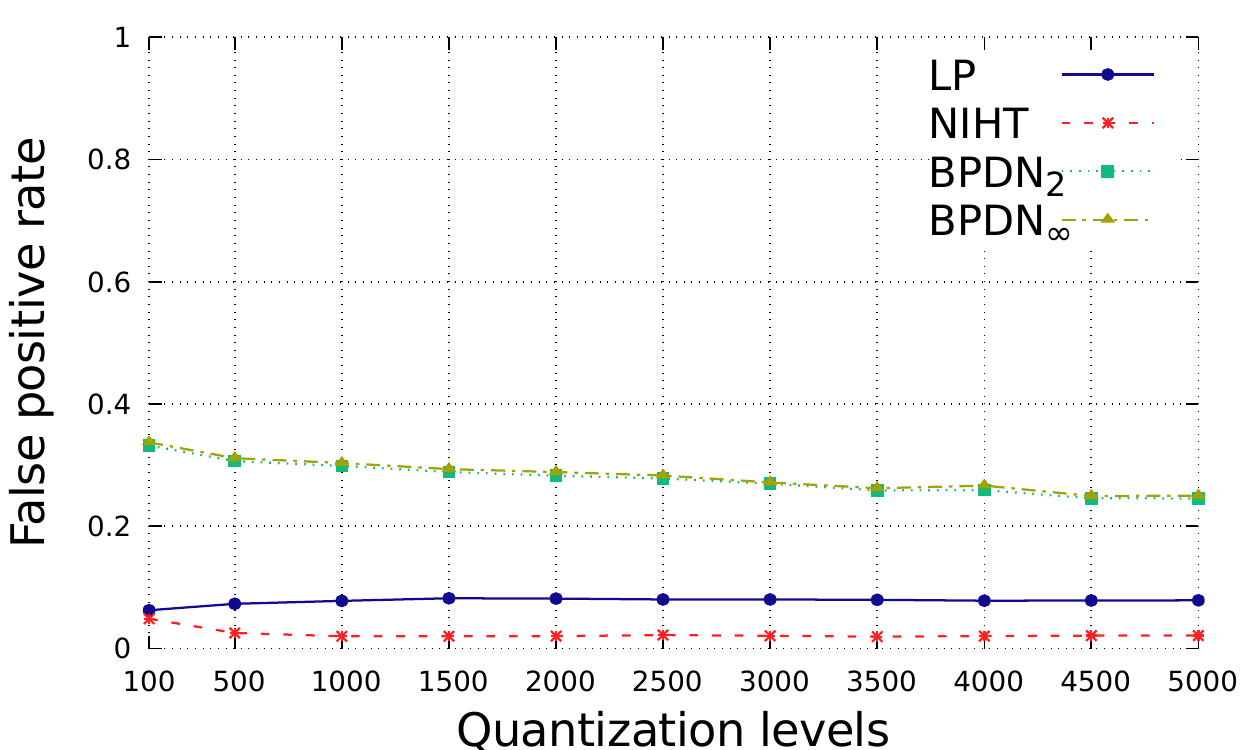} 
\includegraphics[width=0.47\textwidth]{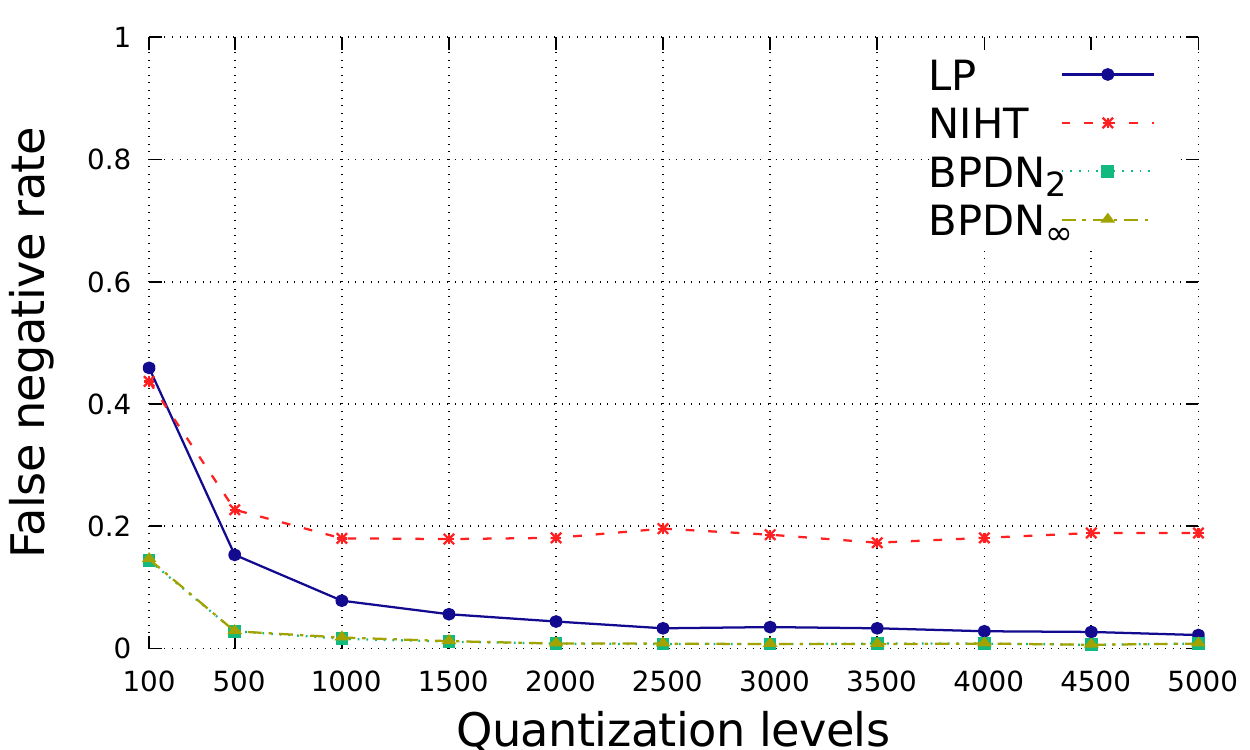}
\caption{Comparison between the proposed LP method and BPDN$_{\infty}$ \cite{don06infty} with $\|\qA x-\qy \|_{\infty}\leq \Dy$, BPDN$_2$ \cite{fou13} with $\|\qA x-\qy \|_{2}\leq \sqrt{m}\Dy$, and NIHT \cite{gur18}.}\label{fig1:errors}
\end{figure*}
\begin{figure*}
\centering
\includegraphics[width=0.47\textwidth]{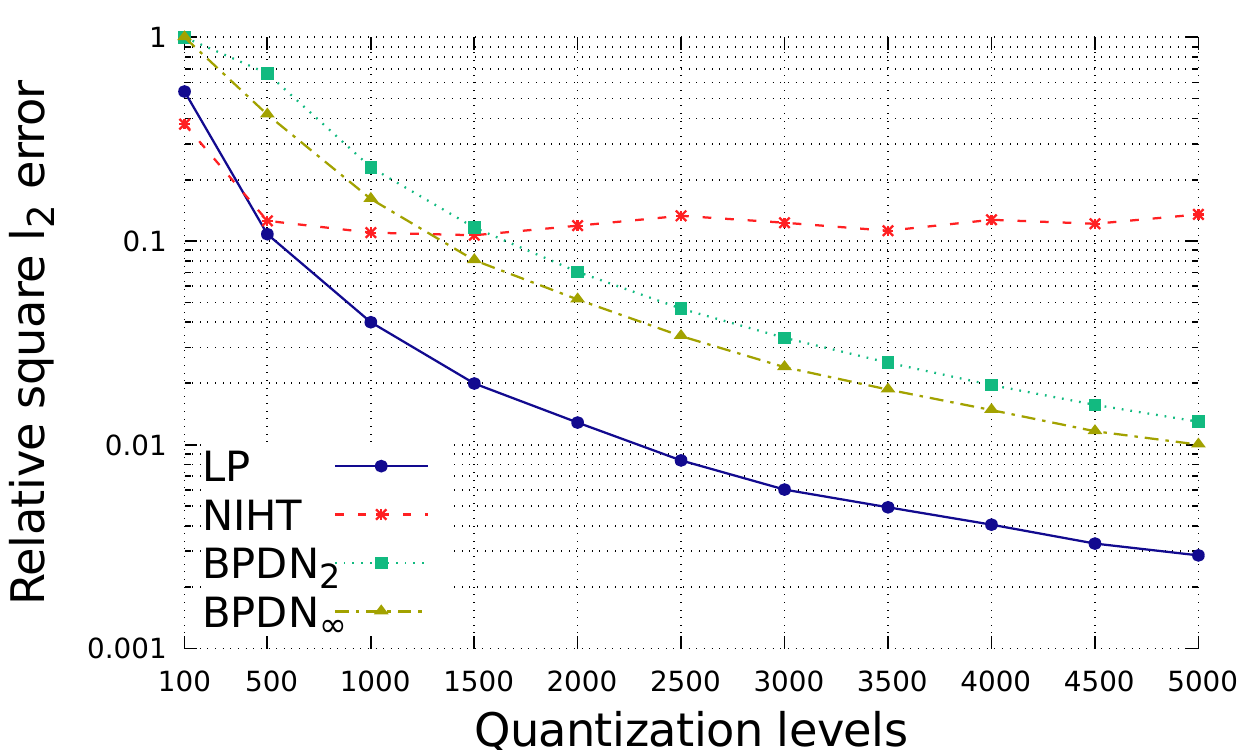}\quad
\includegraphics[width=0.47\textwidth]{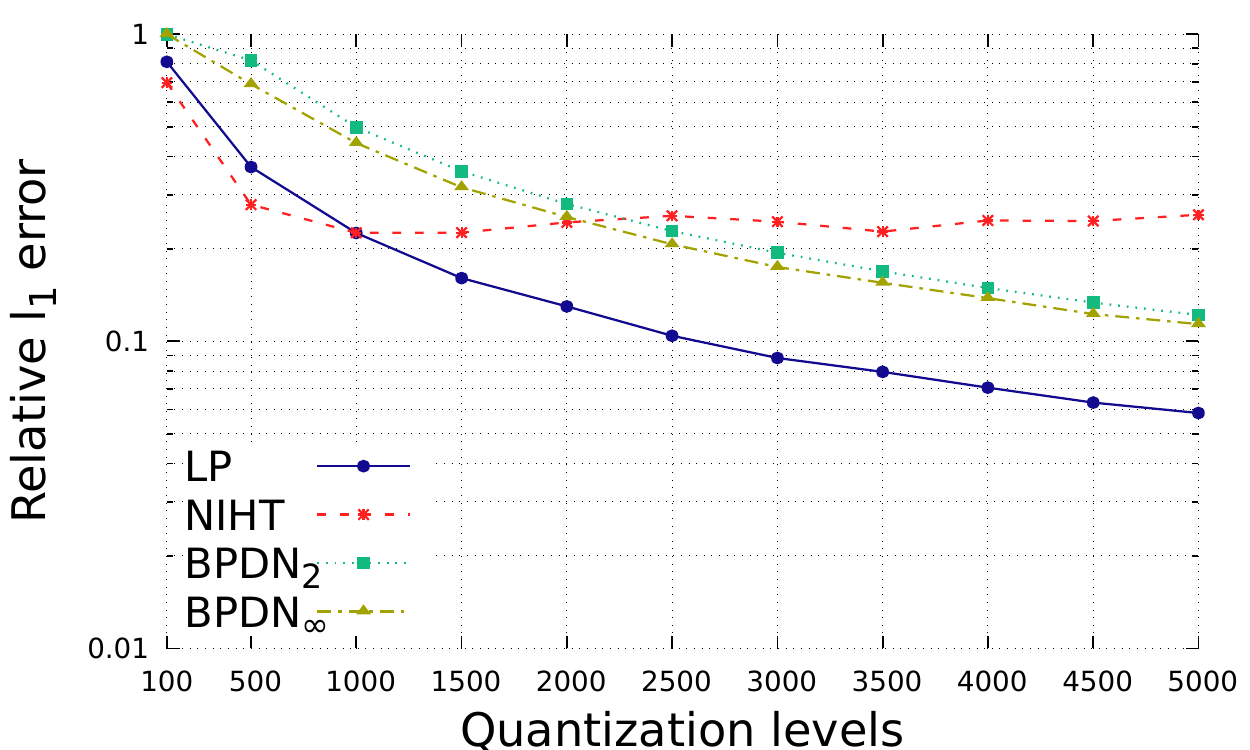}
\includegraphics[width=0.47\textwidth]{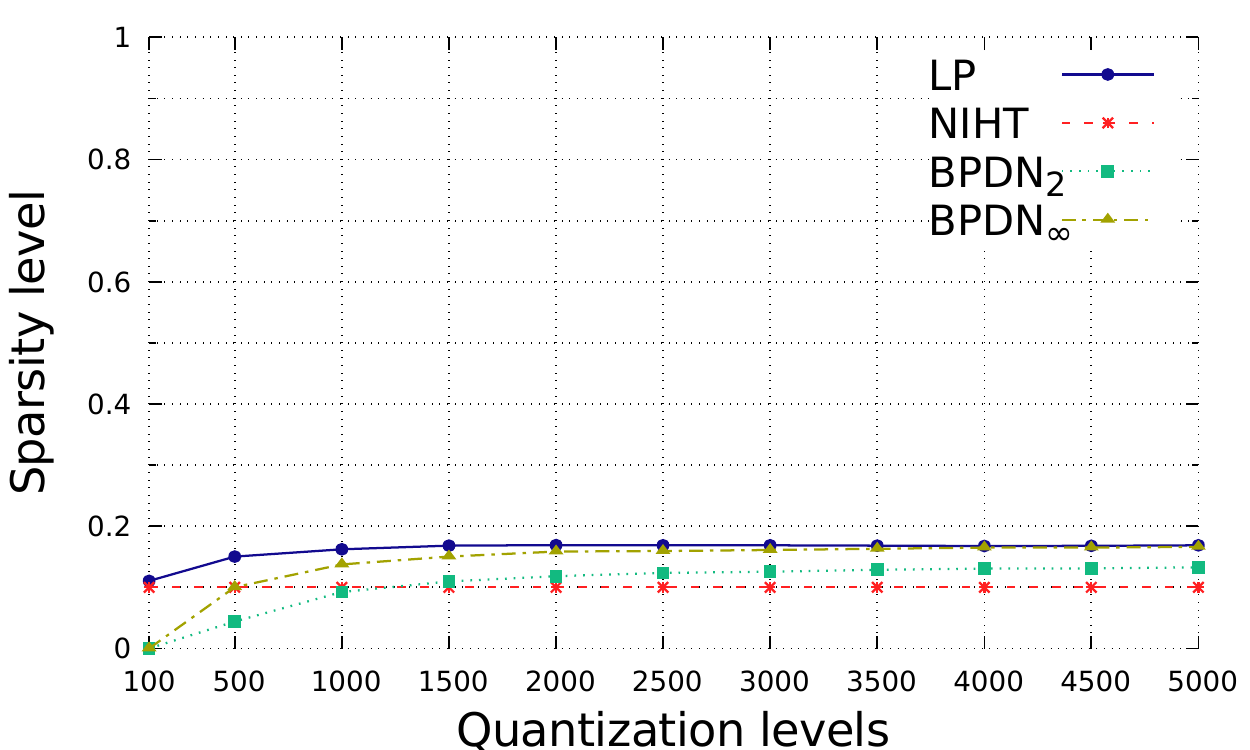}
\includegraphics[width=0.47\textwidth]{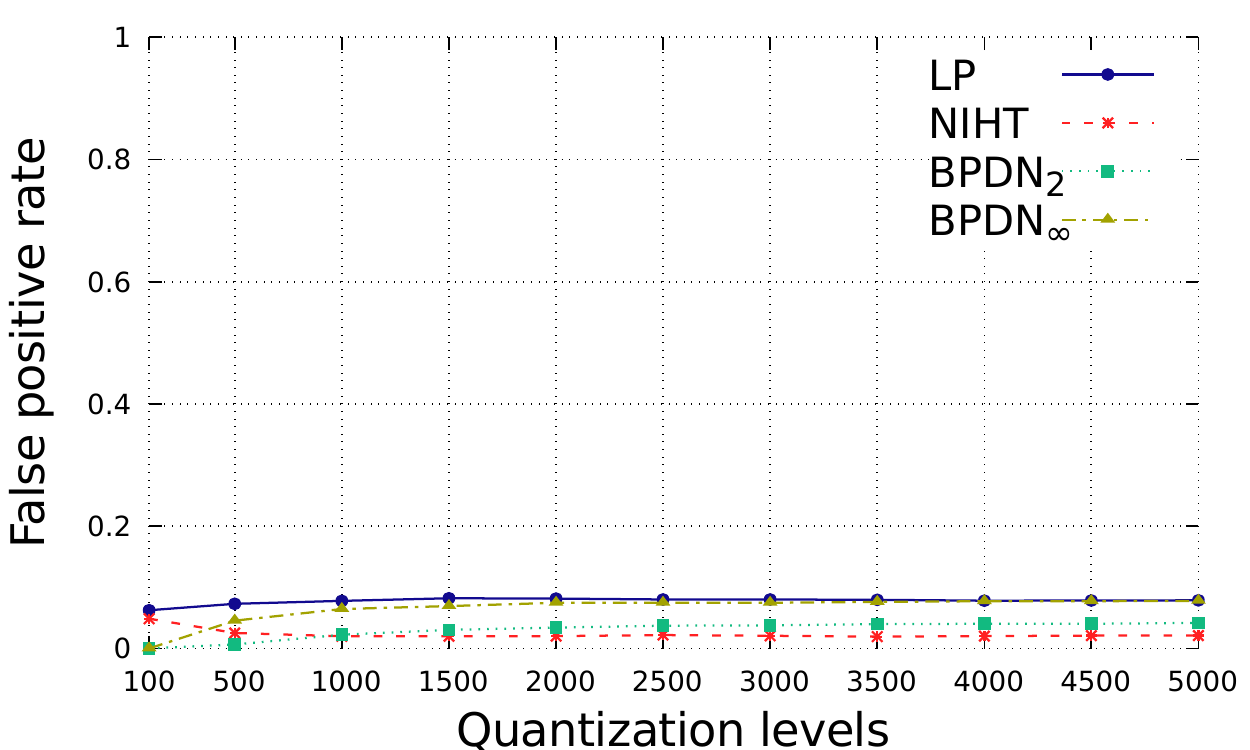} 
\includegraphics[width=0.47\textwidth]{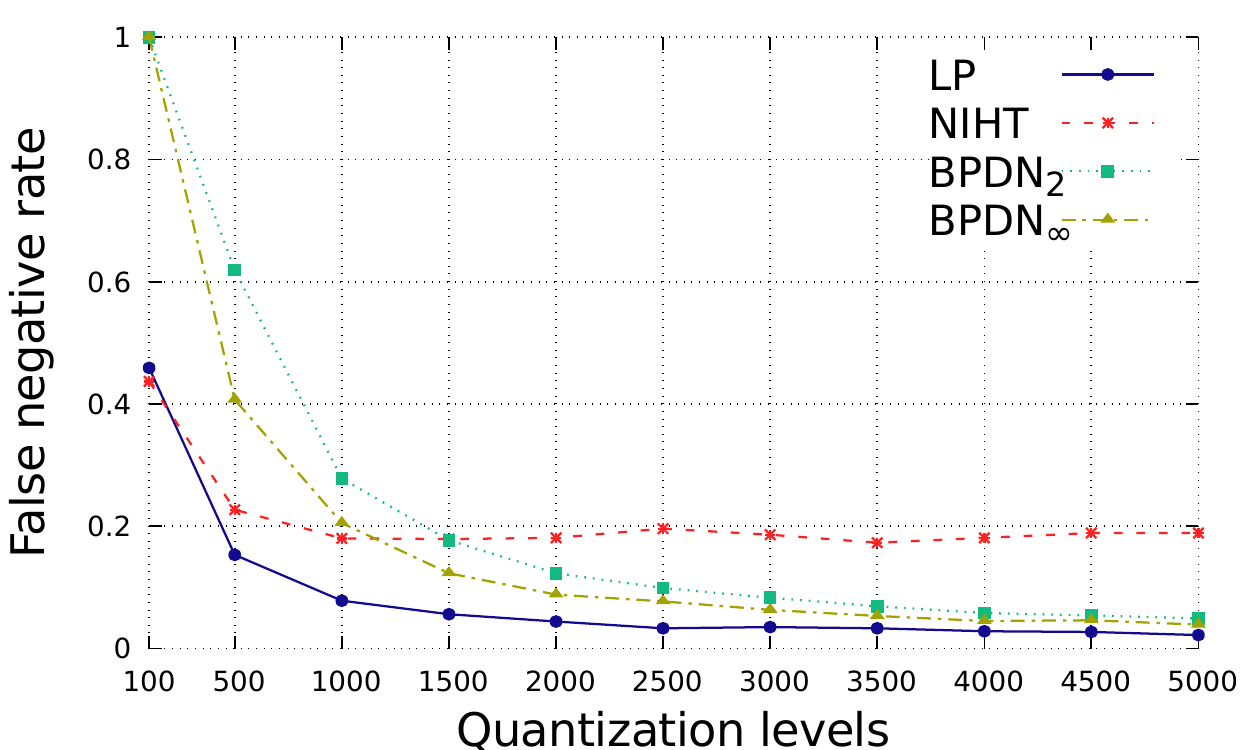}
\caption{Comparison between the proposed LP method and BPDN$_{\infty}$ \cite{don06infty} with $\|\qA x-\qy \|_{\infty}\leq \DA k r+\Dy$, BPDN$_2$ \cite{fou13} with $\|\qA x-\qy \|_{2}\leq \sqrt{m}(\DA k r+\Dy)$, and NIHT \cite{gur18}.}\label{fig2:errors}
\end{figure*}
We propose the following experiment. A system acquires a sparse vector $x\in\R^n$ through a predictor matrix $A\in\R^{m,n}$. The dimensions are $n=100$, $m=40$, $k=10$. $A$ is generated according to a Gaussian distribution $\mathcal{N}(0,\frac{1}{m})$; the support of $x$ is generated uniformly at random, while the non-zero entries uniformly distributed in $[0,r]$ with $r=10$. The quantization is uniform: fixed a certain number of equidistant quantization levels, each entry of $A$ and of $y=Ax$ is approximated with the closest point in the quantization codebook. We assume a quantization range sufficiently large so that saturation problems are negligible. 

We compare the proposed LP method to three state-of-the-art methods:  BPDN$_{\infty}$, BPDN$_{2}$, and the normalized iterative hard thresholding (NIHT),  as presented in \cite{gur18}.

To design the measurement noise bounds for BPDN$_{\infty}$ and BPDN$_2$, we propose two different settings.

- Setting 1: quantization of $A$ is ignored, while we know $\Dy$. Therefore, we impose $\|Ax-y\|_{\infty}\leq \Dy$ and, as a consequence, $\|Ax-y\|_{2}\leq \sqrt{m}\Dy$.

- Setting 2: quantization of $A$ is known, though the related error is moved on the measurements. Assuming to know $\DA$, $k$, and a bound $r>0$ such that $x_i<r$, for any $i=1,\dots,n$, from $\qA x-\qy =\dA x-\dy$ we obtain
$\|\qA x-\qy \|_{\infty}\leq \DA k r+\Dy$
and, as a consequence, 
$\|\qA x-\qy \|_{2}\leq \sqrt{m}(\DA k r+\Dy).$

We notice that NIHT requires the knowledge of $k$.

In Fig. \ref{fig1:errors} and Fig. \ref{fig2:errors}, we show the performance with respect to different quantization levels, from 100 to 5000. For simplicity, we consider that same quantization levels for $A$ and for $y$, while distinguished approximations could be suitably designed. 
The total range is assumed to be $[-r,r]$ with $r=10$. which is generally sufficient to avoid problems of saturation for the proposed setting. We show the relative square $\ell_2$ error, defined as $\|\widehat{x}-x\|_2^2/\|x\|_2^2$; the relative square $\ell_1$ error, defined as $\|\widehat{x}-x\|_1/\|x\|_1$; the normalized sparsity level of the estimation,  the false positive rate, defined as the number of events where $\widehat{x}_i\neq 0$ while $x_i= 0$, over $n-k$; the false negative rate,  defined as the number of events where $\widehat{x}_i= 0$ while $x_i\neq 0$, over $k$.

In Fig. \ref{fig1:errors}, we show the simulations in Setting 1. In this case, the distance from the desired vector is smaller for BDPN$_p$, $p=2,\infty$, with respect to our LP: the smaller feasible set forces a closer consistency to data. However, the smaller feasible set limits the sparsity of the obtained solution, which is between $30\%$ and $40\%$ for BDPN$_p$,  while the correct one is $10\%$. The false positive rate is then large. This is in contrast to the wish of producing sparse solutions.

In order to obtain sparser solutions for BDPN$_p$, $p=2,\infty$, we implement the Setting 2, which is depicted in Fig. \ref{fig2:errors}. By assuming to know $\DA,k,r$, we can suitably enlarge the feasible set and obtain sparse solutions. However, in this case BDPN$_p$, $p=2,\infty$, has a larger number of false negatives, and the recovery accuracy is worse than the proposed LP.


Finally, we notice that the low precision NIHT \cite[Algorithm 1]{gur18} does not show good performance in this experiment. NIHT is proved to be robust when a specific stochastic quantizer is used \cite[Section 3.1]{gur18}, while this experiment shows that further adjustment should be done for other quantizers. We specify that our approach does not require a specific quantization operator and uses only information about the maximum quantization error.

In conclusion, the proposed LP method provides the best performance in this experiment, since at similar sparsity levels, it obtains the best recovery accuracy.
\section{Conclusions}\label{sec:con}
In this paper, we have addressed the problem of sparse linear regression, with particular attention to the compressed case, when only quantized versions of the predictors and of the measurements are available. This problem is relevant in the applications, while  difficult to tackle due to its intrinsic non-convexity. In this work, we have undertaken an $\ell_{\infty}$ approach, based on results in error-in-variables system identification, which allows us to recast the problem into a linear programming model, under suitable sign assumptions. The proposed approach is theoretically proved to be robust, \emph{i.e.}, a finer quantization leads to a smaller recovery error. Moreover, numerical simulations show an improved recovery accuracy with respect to known methods. Generalizations of the sign assumption are possible, as well as the addition of other sources of noise. This will be the subject of future work.
\bibliographystyle{IEEEtran}
\bibliography{refs}
\end{document}